\documentclass{article}                            
\usepackage{amsmath,amssymb,citesort,enumerate,eucal,theorem}

\numberwithin{equation}{section}

\theoremheaderfont{\normalfont\scshape}
{\theorembodyfont{\rmfamily}\newtheorem{dfn}{Definition}[section]}
{\theorembodyfont{\rmfamily}\newtheorem{rmk}{Remark}[section]}
\newtheorem{prop}{Proposition}[section]
\newtheorem{thm}{Theorem}[section]

\DeclareMathOperator{\Bun}{\mathsf{Bun}}
\DeclareMathOperator{\ch}{ch}
\DeclareMathOperator{\Coker}{Coker}
\DeclareMathOperator{\GL}{\mathsf{GL}}
\DeclareMathOperator{\gl}{\mathfrak{gl}}
\DeclareMathOperator{\gr}{gr}
\DeclareMathOperator{\calEnd}{\mathcal{E}nd}
\DeclareMathOperator{\ind}{ind}
\DeclareMathOperator{\Ker}{Ker}
\DeclareMathOperator{\ord}{ord}
\DeclareMathOperator{\Spec}{Spec}
\DeclareMathOperator{\tr}{tr}
\DeclareMathOperator{\virtdim}{virt.dim.}

\def\MgnXb{\mathcal{\overline{M}}_{g,n}(X,\beta)}

\def\bn{\bigskip\noindent}
\def\mn{\medskip\noindent}

\title{Quantum Langlands duality and mirror symmetry}
\author{Igor Yu. Potemine}
\date{}

\begin{document}
\maketitle

\section{Introduction}

There are recent indications that geometric Langlands 
(Beilinson-Drinfeld) duality  \cite{BD:qH,BD:qL} is related to  
$S$-duality in QFT and string theories (\emph{cf.}~\cite{VW}). In this paper 
we discuss a peculiar relation of this Beilinson-Drinfeld duality and 
Givental's version of mirror symmetry. It should be noted that 
mirror symmetry is a $T$-duality from the physical viewpoint \cite{SYZ}.

\medskip Consider two (infinite-dimensional) ind-varieties:
the \emph{$n$-reduced Sato Grassmannian} $\mathcal{GR}^{(n)}$ and the 
\emph{periodic} (or \emph{affine}) \emph{flag variety} $\mathcal{FL}^{(n)}$.
The \emph{Sato Grassmannian} $\mathcal{GR}$ is the set of subspaces 
$W\subset \mathbb{C}(\!(t)\!)$ commensurable with $\mathbb{C}[t^{-1}]$,
\emph{i.~e.}, such that the projection
\begin{equation}
\gamma_W^{}:W\rightarrow \mathbb{C}(\!(t)\!)/t\mathbb{C}[\![t]\!]
\end{equation}
is Fredholm (with finite kernel and cokernel). Then
\begin{equation}
\mathcal{GR}^{(n)}=\big\{W\in\mathcal{GR}\ \big|\ t^{-n}W\subset W\big\}
\end{equation}
and 
\begin{multline}
\mathcal{FL}^{(n)}=\big\{W_0^{}\subset W_1^{}\subset\cdots\subset
W_n\ \big|\ \forall i\ W_i^{}\in\mathcal{GR}^{(n)},\\ \virtdim{W_i}=
i-n\textrm{ and } t^{-n}W_n=W_0\big\}
\end{multline}
where
\begin{equation}
\virtdim{W}=\ind{\gamma_W^{}}=\dim\Ker{\gamma_W^{}}-
\dim\Coker{\gamma_W^{}}
\end{equation}
is the \emph{virtual dimension} of $W$.

\medskip We construct $\mathcal D$-modules of Beilinson-Drinfeld type
on $\mathcal{GR}^{(n)}$ parametrized by local affine $\GL_n$-opers. 
In the same manner, we construct $\mathcal D$-modules of Givental type 
on $\mathcal{FL}^{(n)}$ parametrized by local affine Miura $\GL_n$-opers.

\mn \textsc{Main theorem.} \emph{Givental (quantum) $\mathcal{D}$-modules 
on $\mathcal{FL}^{(n)}$ are transformed to Beilinson-Drinfeld (quantum) 
$\mathcal D$-modules on $\mathcal{GR}^{(n)}$ via the affine Miura 
transformation.}

\medskip The outline of this paper is as follows. First of all, we recall
some facts about the gravitational quantum cohomology and Givental's
construction of a quantum $\mathcal{D}$-module $\mathcal{D}(X)$ for a smooth
projective variety $X$. We would like to construct a family of
$\mathcal{D}$-modules on infinite-dimensional varieties $\mathcal{GR}^{(n)}$ 
and $\mathcal{FL}^{(n)}$ parametrized by affine opers. It is possible to 
generalize Givental's approach using Guest-Otofuji results \cite{GO}.
However, in this paper we adopt a different strategy. In section 4
we recall Beilinson-Drinfeld construction of $\mathcal{D}$-modules 
on the moduli stack $\Bun_G^{}(X)$ of principal $G$-bundles on a curve $X$.
This construction uses a quantization of the Hitchin system. In the 
following two sections (5 and 6) we sketch a local affine version of the 
Beilinson-Drinfeld construction using the formalism of null vectors. 
In some sense, it is related to the quantization of Korteweg-de Vries 
hierarchies. A similar approach is used later to construct 
$\mathcal{D}$-modules on $\mathcal{FL}^{(n)}$ parametrized by local affine 
Miura $\GL_n$-opers. In section 9 we deduce our main theorem from standard 
results about the Miura transformation. Finally, we indicate how one can 
follow very closely the Beilinson-Drinfeld construction using analogues of 
the Hitchin fibration.

\section{Background from gravitational quantum cohomology}

The surveys \cite{CK} and \cite{Ma} may serve as a nice introduction.

\smallskip Physically, (gravitational) quantum cohomology is a 
cohomological realization of non-linear $\Sigma$-models of two-dimensional
field theory (coupled to gravity). One considers interacting ``strings'' 
propagating inside a variety $X$, \emph{i.~e.}, maps $\Sigma\rightarrow X$ 
where $\Sigma$ is a Riemann surface swept out by these strings. The main 
statistical objects of interest are the \emph{correlators} (vacuum 
expectation values) of the \emph{primary fields}, represented by 
cocycles $\gamma_1^{},\cdots, \gamma_n\in\mathrm{H}^*(X,\mathbb Q)$, 
and their \emph{gravitational descendants} 
$\tau_{d_1}^{}\gamma_1^{},\dots,\tau_{d_n}^{}\gamma_n^{}$ 
($d_i\in \mathbb N$). Fix a homology class 
$\beta\in \mathrm{H}_2(X,\mathbb Z)$ and consider the moduli stack
$\MgnXb$ of stable maps $f:C\rightarrow X$ of algebraic curves of genus 
$g$ with $n$ marked points such that $f_*[C]=\beta$. Denote by
$[\MgnXb]^{\mathrm{virt}}$ the \emph{virtual fundamental class} 
(the orbifoldic Poincar\'e dual of 
$1\in\mathrm{H}^*\big(\MgnXb,\mathbb Q\big)$). In addition, let 
\begin{equation}
e_i:\MgnXb\longrightarrow X,\quad e_i^{}(f:C\rightarrow X,p_1,\dots,p_n)=f(p_i)
\end{equation}
be the evaluation maps at $i$th marked point and $\mathcal{L}_i$ 
``the cotangent line at $p_i^{}$'', \emph{i.~e.}, a line bundle over $\MgnXb$ 
whose fiber over $(f:C\rightarrow X,p_1,\dots,p_n)$ is the cotangent
space $T_{p_i}^*C$. Now the \emph{gravitational correlators} are defined
by
\begin{align*}
\langle \tau_{d_1}\gamma_1,\dots,\tau_{d_n}\gamma_n\rangle_{g,\beta}
&=\int_{\MgnXb} I_{g,n,\beta}(\tau_{d_1}\gamma_1,\dots,\tau_{d_n}\gamma_n)\\
\intertext{where}
I_{g,n,\beta}(\tau_{d_1}\gamma_1,\dots,\tau_{d_n}\gamma_n)=&
\big(c_1^{}(\mathcal{L}_1)^{d_1}\cup e_1^*(\gamma_1)\big)\cup\cdots\cup
\big(c_1^{}(\mathcal{L}_n)^{d_n}\cup e_n^*(\gamma_n)\big)
\end{align*}
are the \emph{gravitational Gromov-Witten classes}. If $d_1=\cdots =d_n=0$
then
\begin{align}
I_{g,n,\beta}(\tau_{d_1}\gamma_1,\dots,\tau_{d_n}\gamma_n)
&=I_{g,n,\beta}(\gamma_1,\dots,\gamma_n)\\
\intertext{are the usual \emph{Gromov-Witten classes} and}
\langle \tau_{d_1}\gamma_1,\dots,\tau_{d_n}\gamma_n\rangle_{g,\beta}
=\langle I_{g,n,\beta}\rangle(\gamma_1,\dots,\gamma_n)
\end{align}
are the \emph{Gromov-Witten invariants}.

\section{Givental quantum $\mathcal D$-modules}

First of all, we recall some basic facts about Givental quantum connection
and associated quantum $\mathcal D$-modules. Let $X$ be a smooth projective
variety and $T_0,\dots T_m$ a basis of ${\rm H}^*(X,{\mathbb C})$. We consider
${\rm H}^*(X,{\mathbb C})$ as a Frobenius supermanifold with respect to 
supercommutative variables $t_0^{},\dots,t_m^{}$ associated to $T_i$
(\emph{cf.}~\cite[8.2.1]{CK}). The \emph{Givental quantum connection} 
$\nabla^{\hbar}$ is defined on the trivial cohomology bundle 
${\rm H}^*\big(X,{\mathbb C}[\![t_0,\dots t_m]\!]\big)$ over 
${\rm H}^*(X,{\mathbb C})$ by
\begin{equation}
\nabla_{\frac{\partial}{\partial t_i}}^{\hbar}
\Big(\sum_j a_j^{}T_j^{}\Big)=\hbar\sum_j
\frac{\partial a_j^{}}{\partial t_i^{}} T_j-
\sum_j a_j^{}T_j^{}\ast T_i^{}
\end{equation}  
where $\ast$ is the {\it big quantum product} on 
${\rm H}^*\big(X,{\mathbb C}[\![t_0,\dots t_m]\!]\big)$.

\medskip Consider the \emph{genus $g$ couplings}
\begin{equation}
\langle\langle \tau_{d_1}^{}\gamma_1,\dots,\tau_{d_n}^{}\gamma_n
\rangle\rangle_g=\sum_{k=0}^{\infty}\sum_{\beta} \frac{1}{k!}
\langle\tau_{d_1}^{}\gamma_1,\dots,\tau_{d_n}^{}\gamma_n,
\underbrace{\gamma,\dots,\gamma}_{k\mathrm{ times}}
\rangle_{g,\beta}^{}q^\beta
\end{equation}
where $\gamma=\sum_i t_iT_i$ \cite[sect.~10.1.1]{CK}. 
The formal sections
\begin{equation}
s_a=T_a+\sum_j\left\langle\!\!\!\left\langle\frac{T_a}{\hbar-c},T_j
\right\rangle\!\!\!\right\rangle_0T^j,
\end{equation}
where $c=c_1^{}(\mathcal{L}_1)$, form a basis of $\nabla^{\hbar}$-flat
sections of ${\rm H}^*\big(X,{\mathbb C}[\![t_0,\dots t_m]\!]\big)$.
This means that
\begin{equation}
\hbar\frac{\partial s_a}{\partial t_i}=T_i\ast s_a,\quad
a,i=0,\dots,m
\end{equation}
\cite[sect.~10.2.1]{CK}. 

\begin{dfn} Denote
\begin{equation}
\mathcal D=\mathbb C\left[\hbar \partial /\partial t_i,\exp t_i,\hbar
\right].
\end{equation}
The elements of $\mathcal{D}$, called \emph{quantum differential operators},
 are polynomials in the quantities
\begin{equation}
\hbar \partial /\partial t_0,\dots,\hbar \partial /\partial t_r,
e^{t_0},\dots,e^{t_r},\hbar.
\end{equation}
The $\mathcal{D}$-module 
\begin{equation}
\mathcal D(X)=\mathcal D/I,\quad I=\{D\in\mathcal D\;|\;
DS_a=0,\ 1\leqslant a\leqslant r \} 
\end{equation}
generated by $S_a=\langle s_a,1\rangle$ is called the 
\emph{Givental quantum $\mathcal D$-module of $X$}.
\end{dfn}

\medskip Givental and Kim showed that the quantum $\mathcal{D}$-module
of the flag variety $G/B$ is generated by non-constant first integrals 
of the quantum Toda chain for the Langlands dual $G^{\vee}$
\cite{Gi,Kim}. Guest and Otofuji \cite{GO} generalized some of these
results to the case of the affine flag variety $\mathcal{FL}^{(n)}$.

\medskip We use this cohomological Givental construction as a motivation and
adopt a different strategy in the rest of the paper. 

\section{Hitchin systems, $G$-opers and Beilinson-Drin\-feld duality}

Let $X$ be a smooth projective curve of genus $g\geqslant 2$ and 
$\Bun_G^{}(X)$ the moduli space of stable principal $G$-bundles 
on $X$. If $G=\GL_n$ we denote $\Bun_{X,n}^{}$ the moduli space 
of stable vector bundles of rank $n$. Hitchin \cite{Hi} constructed
a map  
\begin{equation}
h:T^*\Bun_G^{}(X)\rightarrow V
\end{equation}
from the cotangent bundle $T^*\Bun_G^{}(X)$ to a certain
vector space $V$ of dimension
\begin{equation}
\dim V=\dim\Bun_G^{}(X)=(g-1)\dim G+\dim Z(G).
\end{equation}

\smallskip If $G=\GL_n$ then
\begin{equation}
\big(\Bun_G^{}(X)\big)_P\simeq\mathrm{H}^0(X,\calEnd P\otimes K_X)
\end{equation}
where $K_X$ denotes the canonical bundle of $X$. In other words, 
$T^*\Bun_{X,n}^{}$ consists of pairs $(P,\varphi)$ where 
$P\in \Bun_{X,n}^{}$ and $\varphi:P\rightarrow P\otimes K_X$
a twisted endomorphism. The Hitchin map is defined by
\begin{align}
h:T^*\Bun_{X,n}^{}&\longrightarrow V=\bigoplus_{i=1}^n 
\mathrm{H}^0(X,K_X^i)\\
(P,\varphi)&\longmapsto \ch(\varphi)=\big(\tr\varphi,\tr(\wedge^2\varphi),
\dots,\tr(\wedge^n\varphi)\big)
\end{align}
where $\ch(\varphi)$ is the $n$-tuple of the coefficients of the
characterstic polynomial of~$\varphi$. 

\smallskip Hitchin has shown that for any basis $\{f_1,\dots,f_r\}$ of 
$V^*$ the functions $h^*f_i$ and $h^*f_j$
commute with respect to the standard Poisson bracket: 
\begin{equation}
\{h^*f_i,h^*f_j\}=0
\end{equation}
on $T^*\Bun_G^{}(X)$. Denote by $B_G(X)$ the ring of polynomial functions 
on $V$. The Hitchin map induces a morphism
\begin{equation}
\pi_{\mathrm{cl}}^{}:B_G(X)\rightarrow\{\textrm{functions on }
T^*\Bun_G^{}(X)\}
\end{equation}
Suppose now, for the sake of simplicity, that $G$ is semi-simple 
and simply connected. Beilinson and Drinfeld constructed a 
\emph{quantization} (\emph{$\hbar$-deformation}) of $\pi_{\mathrm{cl}}$, 
that is, a graded $\mathbb C$-algebra $A_G(X)$ such that $\gr A_G(X)=B_G(X)$
equipped with a morphism
\begin{equation}
\pi:A_G(X)\longrightarrow\mathrm{H}^0\big(\Bun_G^{}(X),\mathcal{D}^{1/2}
\big)=\mathrm{H}^0\big(\Bun_G^{}(X),\mathcal{D}_{-h^{\vee}}\big)
\end{equation}
compatible with $\pi_{\mathrm{cl}}^{}$. Here $h^{\vee}$ is the dual 
Coxeter number of $G$ and $\mathcal{D}^{1/2}=\mathcal{D}_{-h^{\vee}}$
is the sheaf of twisted differential operators (of critical level
$k=-h^{\vee}$) acting on $K_{\Bun_G^{}(X)}^{1/2}=\eta^{-h^{\vee}}$,
where $\eta^k$, $k\in \mathbb Z$, denotes a line bundle on $\Bun_G^{}(X)$
obtained via the generalized Borel-Weil-Bott (Kumar-Mathieu) theorem
\cite{BD:qH,BD:qL,Fr:aff}.


\medskip The elements of $A_G(X)$ are called \emph{$G^{\vee}$-opers}.
Denote by $\mathcal{M}\big(\Bun_G^{}(X)\big)$ (resp. 
$\mathcal{M}\big(X\times\Bun_G^{}(X)\big)$ the category of
$\mathcal D$-modules on $\Bun_G^{}(X)$ (resp. $X\times\Bun_G^{}(X)$).

\begin{thm}[Beilinson-Drinfeld]
There\hskip 2mm exists\hskip 2mm a\hskip 2mm family\hskip 2mm
of\hskip 3mm $\mathcal D$-modules\\ $\{\mathcal{D}_{\mathcal L}\}$ 
on $\Bun_G^{}(X)$ parametrized by $G^{\vee}$-opers 
$\{\mathcal{L}\}$ which are eigensheaves of Hecke operators:
\begin{equation}
T_{\chi}:\mathcal{M}\big(\Bun_G^{}(X)\big)\longrightarrow 
\mathcal{M}\big(X\times \Bun_G^{}(X)\big),\quad T_{\chi}
\mathcal{D}_{\mathcal L}=V_{\mathcal L}^{\chi}\boxtimes
\mathcal{D}_{\mathcal L}
\end{equation}
for any dominant coweight $\chi\in P_+(G^{\vee})$ where 
$V_{\mathcal L}^{\chi}$ is a $G^{\vee}$-local system on $X$
associated to a couple $(\mathcal{L},\chi)$.
\end{thm}

\medskip One can give an explicit description of $\GL_n$-opers.

\begin{dfn} A \emph{$\GL_n$-oper} over $X$ is a local system 
$(\mathcal{E},\nabla)$ adapted to a complete flag
\begin{equation}
0=\mathcal{E}_0\subset \mathcal{E}_1\subset\cdots\subset \mathcal{E}_n=
\mathcal{E}
\end{equation} 
of vector bundles over $X$ such that
\begin{enumerate}[(i)]
\item $\gr_i^{}(\mathcal{E})$ are invertible,
\item $\nabla(\mathcal{E}_i)\subset \mathcal{E}_{i+1}\otimes K_X$,
\item $\nabla$ induces an isomorphism
\begin{equation}
\gr_i^{}(\mathcal{E})\overset{\sim}{\longrightarrow}\gr_{i+1}^{}(\mathcal{E})
\otimes K_X.
\end{equation}
\end{enumerate}
\end{dfn}

In the local coordinates a $\GL_n$-oper is of the form:
\begin{equation}
\partial_t+
\begin{pmatrix}
* &* &\cdots &* \cr
+ &* &\ddots &\vdots \cr
&\ddots &\ddots & ^{\textstyle *} \cr
0 &&+ &*
\end{pmatrix}
\end{equation}
where elements denoted by $*$ are arbitrary and elements denoted by $+$ are 
non-zero. Locally, one can write a $\GL_n$-oper in the form
\begin{equation}\label{eq:locoper}
\partial_t-
\begin{pmatrix}
q_1^{} &q_2^{} &\cdots &q_n^{} \cr
1 &0 &\cdots &0 \cr
&\ddots &\ddots &\vdots \cr
0 &&1 &0
\end{pmatrix},
\; q_i^{}\in\mathbb{C}[\![t]\!],
\end{equation}
equivalent to a differential operator
\begin{equation}
L=\partial_t^n-q_1^{}\partial_t^{n-1}-\cdots-q_n^{}.
\end{equation}
Thus, a \emph{local $\GL_n$-oper} is a differential operator of order $n$ 
whose principal symbol is equal to 1.

\section{Affine $\GL_n$-opers and KdV hierarchies}\label{sect:opers}

We suppose now that $X$ is a complete (possibly singular) curve and 
$P$ is a smooth closed point on $X$.

\begin{dfn} A \emph{local affine $\GL_n$-oper} is a pair 
$\big(\widehat{\mathcal{E}},\widehat{\nabla}\big)$ where 
$\widehat{\mathcal{E}}$ is a formal torsion-free sheaf of rank $n$ 
over $X\times\widehat{D}_t=X\times\Spec \mathbb{C}[\![t]\!]$
(an infinitesimal deformation of a torsion-free sheaf $\mathcal{E}$ 
of rank $n$ over $X$) equipped with a parabolic structure along of 
$P\times \widehat{D}_t$:
\begin{equation}
\widehat{\mathcal{E}}=\widehat{\mathcal{E}}_0\subset
\widehat{\mathcal{E}}_1\subset\cdots\subset
\widehat{\mathcal{E}}_n=\widehat{\mathcal{E}}(P)
\end{equation}
Moreover, $\widehat{\nabla}$ is a (micro)connexion on $\widehat{\mathcal{E}}$, 
that is, a map
\begin{equation}
\widehat{\nabla}:\left(\varinjlim \widehat{\mathcal{E}}_i\right)
\longrightarrow \left(\varinjlim \widehat{\mathcal{E}}_i\right)
\end{equation}
satisfying the Leibniz rule
\begin{equation}
\widehat{\nabla}(fs)=f\widehat{\nabla}(s)+
\frac{\partial f}{\partial t}\cdot s
\end{equation}
for any function $f$ on $X\times\widehat{D}_t$ and any section $s$,
and such that

\begin{enumerate}[(i)]
\item $\gr_i^{}\big(\widehat{\mathcal{E}}\big)$ are invertible
and $\chi\big(\widehat{\mathcal{E}}_{i+1}\big)=
\chi\big(\widehat{\mathcal{E}}_i\big)+1$,
\item $\widehat{\nabla}\big(\widehat{\mathcal{E}}_i\big)\subset
\widehat{\mathcal{E}}_{i+1}$
and $\ord_P^{}\big(\widehat{\nabla}s\big)=\ord_P^{}(s)+1$, and
\item $\gr_i^{}\big(\widehat{\mathcal{E}}\big)\overset{\sim}{\longrightarrow}
\gr_{i+1}^{}\big(\widehat{\mathcal{E}}\big)$ is an isomorphism induced by 
$\widehat{\nabla}$.
\end{enumerate}
\end{dfn}

\begin{prop}[Drinfeld \cite{Dr:comm}] There exists a natural bijection between 
local affine $\GL_n$-opers and \emph{Krichever modules}, \emph{i.~e.},
commutative subrings $R\subset \mathbb{C}[\![t]\!][\partial_t]$,
of rank $n$.
\end{prop}

\begin{rmk} In positive characteristic, this corresponds to a bijection
between elliptic sheaves over a field and Drinfeld modules.
\end{rmk}

The famous Krichever map \cite{Mul} associates to any 
local affine $\GL_n$-oper a point $W$ on the $n$-reduced Sato 
Grassmannian $\mathcal{GR}^{(n)}$. 

\medskip The $n$th Korteweg-de Vries (KdV) hierarchy describes  
isospectral deformations of a differential operator $L\in\mathbb{C}
[\![t]\!][\partial_t]$ of order $n$. It may be written in the Lax form:
\begin{equation}
\frac{\partial L}{\partial t_r}=\big[L_+^{r/n},L\big],\;r\in\mathbb{N},
\end{equation}
where $L_+^{r/n}$ is the positive part of the microdifferential operator
\begin{equation}
L^{r/n}\in\mathbb{C}[\![t]\!]\big(\!\!\big(\partial_t^{-1}\big)\!\!\big).
\end{equation}
Indeed, by a lemma of Schur \cite[p.~140]{Mum}, any differential operator
$L\in\mathbb{C}[\![t]\!][\partial_t]$ of order $n$ has an (essentially
unique) root $L^{1/n}\in\mathbb{C} [\![t]\!]
\big(\!\!\big(\partial_t^{-1}\big)\!\!\big)$.

\smallskip One can view the $n$th KdV hierarchy as a dynamical system on 
the space of local affine $\GL_n$-opers.

\section{Semi-infinite $q$-wedges, quantum $\tau$-functions
and quantum $\mathcal{D}$-modules for KdV theory}\label{sect:qwedges}

Consider the infinite-dimensional Clifford algebra 
$\mathcal{CL}$ generated by $\psi_n^{\pm}$, $n\in
{\mathbb Z}+1/2$, satisfying the following relations:
\begin{equation}
[\psi_m^{\pm},\psi_n^{\pm}]_+=0,\quad
[\psi_m^+,\psi_n^-]_+=\delta_{m,-n}.
\end{equation}
This Clifford algebra admits a unique irreducible 
representation (called \emph{spin representation}) in
an infinite-dimensional vector space ${\mathcal V}$ (resp. ${\mathcal V}^*$) 
which is a left (resp. right) module admitting a non-zero 
\emph{vacuum vector} $|0\rangle$ (resp. $\langle0|$) such that
\begin{equation}
\psi_i^{\pm}|0\rangle=0\text{ (resp. }\langle 
0|\psi_{-i}^{\pm}=0) \textrm{ for } i>0.
\end{equation}
Denote by $\gl_{\infty}$ the Lie algebra of infinite
matrices indexed by ${\mathbb Z}+1/2$ with almost all entries 
equal to zero. It has a canonical basis consisting of
matrices $E_{ij}$, $i,j\in {\mathbb Z}+1/2$, with 1 as the
$(i,j)$-entry and 0 elsewhere. 
The map $E_{ij}\mapsto \psi_{-i}^{+}\psi_j^{-}$ defines a 
bijection between $\gl_{\infty}$ and quadratic elements 
$\sum a_{ij}^{}\psi_{-i}^{+}\psi_j^{-}$ of $\mathcal{CL}$.
Later $\gl_{\infty}$ will also be identified with normally ordered
expressions of the type $\sum a_{ij}^{}{\rm :}\psi_{-i}^{+}\psi_j^{-}{\rm :}$.

\smallskip Using the exponential map we obtain a representation $R$ of the 
group $\GL_\infty$ on ${\mathcal V}$ and ${\mathcal V}^*$. For a positive 
integer $m$ denote
\begin{equation*}
\langle\pm m|=\langle 0|\psi_{1/2}^{\pm}\cdots\psi_{m-1/2}^{\pm}
\in {\mathcal V}^* \textrm{ and } |\pm m\rangle=\psi_{-m+1/2}^{\pm}\cdots
\psi_{-1/2}^{\pm}|0\rangle\in {\mathcal V}.
\end{equation*}
This defines the \emph{charge decomposition} ${\mathcal V}=\oplus 
\mathcal{V}^{(m)}$ (\emph{cf.} below). 

\smallskip This spin representation has the following semi-infinite wedge 
realization. Let $\mathbb{C}^{\infty}$ be an infinite-dimensional complex
vector space with a fixed basis $\{v_i^{}\}_{i\in{\mathbb Z}+1/2}$. Then the 
\emph{semi-infinite wedge} (or \emph{fermionic Fock}) \emph{space} 
${\mathcal V}=\wedge^{\infty/2}\mathbb{C}^{\infty}$ is generated 
by semi-infinite monomials of the form 
\begin{equation}
\begin{split}
v_{i_1^{}}\wedge v_{i_2^{}}\wedge v_{i_3^{}}\wedge\cdots 
&\quad \hbox{where}\  i_1^{}>i_2^{}>i_3^{}>\cdots\\ &\quad\textrm{and } 
i_{k+1}^{}=i_k^{}-1 \textrm{ for } k\gg0.
\end{split}
\end{equation}
The generators of the infinite-dimensional Clifford algebra
are represented as the wedging and contracting operators defined by:
\begin{align*}
&\psi_j^+(v_{i_1}^{}\wedge v_{i_2}^{}\wedge\cdots)
=v_{-j}^{}\wedge v_{i_1}^{}\wedge v_{i_2}^{}\wedge\cdots\\
&\hskip 12mm=(-1)^s v_{i_1^{}}^{}\wedge\cdots\wedge v_{i_s^{}}^{}
\wedge v_{-j}^{}\wedge v_{i_{s+1}^{}}^{}\cdots\quad 
\textrm{if } i_s^{}>-j>i_{s+1}^{}\\
&\psi_j^-(v_{i_1}^{}\wedge v_{i_2}^{}\wedge\cdots)
=(-1)^{s+1}v_{i_1^{}}^{}\wedge\cdots\wedge v_{i_{s-1}^{}}^{}
\wedge v_{i_{s+1}^{}}^{}\cdots\quad \textrm{if } j=i_s^{}
\end{align*}
and zero in all other cases. By definition, the \emph{charge} of 
$v_{i_1^{}}\wedge v_{i_2^{}}^{}\wedge\cdots$ is equal to $m$ if
$i_k^{}+k=m+1/2$ for $k\gg0$. This defines the charge decomposition
as above. In particular,
\begin{equation}
|m\rangle=v_{m-1/2}^{}\wedge v_{m-3/2}^{}\wedge v_{m-5/2}^{}
\cdots\in {\mathcal V}^{(m)}.
\end{equation}
The space ${\mathcal V}^{(m)}$ is an irreducible highest weight representation
of $\gl_{\infty}\subset \mathcal{CL}$ with the highest weight vector
$|m\rangle$. 

\smallskip The $\GL_{\infty}$-orbit of the vacuum vector $|0\rangle$
is closely related to the infinite-dimensional Sato Grassmannian. 
Indeed, let us identify $\mathbb{C}^{\infty}$ with $\mathbb{C}(\!(t)\!)$
via $v_i^{}\leftrightarrow t^{i+1/2}$. Then the map
\begin{equation*}
\varphi:\GL_{\infty}|0\rangle\longrightarrow \mathcal{GR},
\ \tau=\bigwedge_{\{i=-1/2,-3/2,\cdots\}} u_i^{}\longmapsto 
\sum_{\{i=-1/2,-3/2,\cdots\}} \mathbb{C}u_i^{}
\end{equation*}
is well-defined since $u_{-i}^{}=v_{-i}^{}$ for $i$ big enough.
This map is surjective and we can also write
\begin{equation}
\mathcal{GR}=\GL_\infty|0\rangle/\mathbb{C}^*.
\end{equation}
By definition, semi-infinite wedges lying in the
$\GL_{\infty}$-orbit of $|0\rangle$ are called 
\emph{$\tau$-functions} in the fermionic picture.

\smallskip Consider the \emph{affine Hecke} (or \emph{Iwahori-Hecke})
\emph{algebra} $\mathcal{\widehat{H}}_N(q)$ with the generators $T_i^{\pm1}$,
$X_i^{\pm1}$, $1\le i\le N-1$, and relations:
\begin{subequations}
\begin{align}
&T_i^{}T_i^{-1}=T_i^{-1}T_i^{}=X_i^{}X_i^{-1}=X_i^{-1}X_i^{}=1;\\
&T_iT_j=T_jT_i\textrm{ if } |i-j|>1,\quad X_iX_j=X_jX_i\ \forall i,j;\\
&T_iT_{i+1}T_i=T_{i+1}T_iT_{i+1},\quad (T_i+1)(T_i-q)=0;\\
&X_jT_i=T_iX_j \textrm{ if } j\ne i,i+1, \quad T_iX_iT_i=qX_{i+1}^{}. 
\end{align}
\end{subequations}
Its infinite counterpart $\mathcal{\widehat{H}}_{\infty}(q)$ is generated by
a countable number of generators $T_i^{\pm1}$ and $X_i^{\pm1}$, 
$i\in{\mathbb N}$, satisfying the same relations. 

\smallskip Let $V$ be a complex vector space of dimension $n$ with a
basis $\{e_1^{},\dots,e_n^{}\}$ and let $V(\mathfrak z)=
V\otimes_{\mathbb{C}}^{}\mathbb{C}[\mathfrak z,\mathfrak z^{-1}]$. 
We identify $V(\mathfrak z)$ with $\mathbb{C}^{\infty}$ via the correspondence 
$\mathfrak{z}^je_i^{}\leftrightarrow v_{i-nj-1/2}$. There are 
well-known right actions of $\mathcal{\widehat{H}}_N(q^2)$ on 
$V(\mathfrak z)^{\otimes N}$ and of $\mathcal{\widehat{H}}_{\infty}(q^2)$ 
on the so-called \emph{thermodynamic limit} 
$V(\mathfrak z)^{\otimes\infty/2}$ \cite{St,KMS}. 
The latter is well-defined since any $T_i$ acts only in a pair of 
adjacent factors. The $q$-antisymmetrization procedure (loc.$\;$cit.) 
gives the $q$-wedge spaces
\begin{subequations}
\begin{align}
\bigwedge_q^N V(\mathfrak z) &=V(\mathfrak z)^{\otimes n}\left/
\sum_{i=1}^{n-1}\Ker(T_i+1)\right.\\
\intertext{and}
\bigwedge_q^{\infty/2} V(\mathfrak z) &=V(\mathfrak z)^{\otimes\infty/2}
\left/\sum_{i=1}^{\infty} \Ker(T_i+1)\right..
\end{align}
\end{subequations}
In the same manner, the $q$-deformed Fock space ${\mathcal V}_q^{(m)}$ 
of charge $m$ is the quotient of the space $U^{(m)}$ of pure tensors 
$u_{m_1^{}}^{}\otimes u_{m_1^{}}^{}\otimes\cdots$ where $m_k^{}=m-k+1/2$ 
for $k\gg1$ by $\sum_{i=1}^{\infty} (T_i+1)$.

\begin{dfn} The elements of the $q$-antisymmetrized
$\GL_{\infty}$-orbit of the vaccum vector $|0\rangle$ will be 
called \emph{quantum $\tau$-functions}.
\end{dfn}

\medskip Using fermionic normal ordering 
\begin{equation}
{\rm :}\psi_m^{}\psi_n^*{\rm :}=
\begin{cases}
\psi_m^{}\psi_n^* &\text{if $m<0$ or if $n>0$},\\
-\psi_n^*\psi_m &\text{if $m>0$ or if $n<0$}
\end{cases}
\end{equation}
one can define operators 
\begin{equation}
H_n=\sum_{i\in {\mathbb Z}+1/2} {\rm :}\psi_{-i}^+
\psi_{i+n}^-{\rm :}
\end{equation}
satisfying the commutation relations:
\begin{equation}
[H_m,H_n]=m\delta_{m,-n}^{}.
\end{equation}
It is easy to see that for any given element $u$ of the Fock space,
the expression $H_n|u\rangle$ is a finite sum. Consider also the generating
function
\begin{equation}
H(t)=\sum_{n=1}^{\infty} t_n^{}H_n{}
\end{equation}
called the \emph{Hamiltonian} where $t=(t_1,t_2,\dots)$. In the bosonic 
picture, the $\tau$-function $\tau_g$, associated to 
$g=\sum a_{ij}^{}{\rm :}\psi_{-i}^{+}\psi_j^{-}{\rm :}\in\gl_{\infty}$,
is the correlation function
\begin{equation}
\tau_g(t)=\langle\,0\,|\,e^{H(t)}g\,|\,0\,\rangle=\langle\,0\,|\,e^{H(t)}
e^{\sum a_{ij}^{}{\rm :}\psi_{-i}^{+}\psi_j^{-}{\rm :}}\,|\,0\,\rangle
\end{equation}
(\emph{cf.}~\cite{DJKM,DJM,JM}). The wedging and contracting operators 
$\psi_i^+$ and $\psi_j^-$ correspond to the intertwiners: 
\begin{equation}\label{geomtwins1}
\Gamma_i^+(\mathfrak z):V(\mathfrak z)\otimes \mathcal{V}^{(m)}
\longrightarrow \mathcal{V}^{(m+1)}\textrm{ and } \Gamma_j^-(\mathfrak z):
\mathcal{V}^{(m)}\longrightarrow \mathcal{V}^{(m-1)}\otimes 
V(\mathfrak z).
\end{equation}

Now let $W$ be a point of $\mathcal{GR}$ associated to a local affine 
$\GL_n$-oper $\big(\widehat{\mathcal{E}},\widehat{\nabla}\big)$ via the 
Krichever map (\emph{cf.}~sect.~\ref{sect:opers}). Consider
the associated $\tau$-function
\begin{equation}
\tau_W^{}(t)=\langle\,0\,|\,e^{H(t)}g_W^{}\,|\,0\,\rangle=\langle\,0\,
|\,e^{H(t)}e^{\sum a_{ij}^{}{\rm :}\psi_{-i}^{+}\psi_j^{-}{\rm :}}
\,|\,0\,\rangle
\end{equation}
The wedging and contracting operators $\psi_i^+$ and $\psi_j^-$ correspond 
to the intertwiners of level-1 $\widehat{\mathfrak{sl}}_n$-modules
\begin{equation}\label{geomtwins2}
\Gamma_i^+(\mathfrak z):V(\mathfrak z)\otimes\mathcal{V}_{\Lambda_m}
\longrightarrow\mathcal{V}_{\Lambda_{m+1}}\textrm{ and }
\Gamma_j^-(\mathfrak z):\mathcal{V}_{\Lambda_m}\longrightarrow 
\mathcal{V}_{\Lambda_{m-1}}\otimes V(\mathfrak z)
\end{equation}
where $V_{\Lambda_m}$ are the $\widehat{\mathfrak{sl}}_n$-modules
with fundamental weights $\Lambda_m$.

\medskip In general, consider a \emph{singular vector}
of a highest weight $\widehat{\mathfrak{sl}}_n$-module $V_\Lambda$. 
It is an eigenvector of the Cartan subalgebra of $\widehat{\mathfrak{sl}}_n$ 
annihilated by $U_+\big(\widehat{\mathfrak{sl}}_n\big)$ different from the
highest weight vector $v_{\Lambda}^{}$. It may be written in the form 
$Sv_{\Lambda}^{}$ where $S$ is a uniquely defined element of 
$U_-\big(\widehat{\mathfrak{sl}}_n\big)$. Then
\begin{equation}
\langle\,0\,|\,(Sv_\Lambda^{})g\,|\,0\,\rangle=
D\langle\,0\,|\,v_\Lambda^{}g\,|\,0\,\rangle=0.
\end{equation} 
for a certain differential operator $D$ \cite[sect.~5.3]{FM}.

\medskip In our case, any $\tau$-function $\tau_W^{}$ defines a
$\mathcal{D}$-module $\mathcal{D}(\tau_W^{})$ on $\mathcal{GR}$. 
These $\mathcal{D}$-modules are closely related to the so-called 
deformed Knizhnik-Zamolodchikov equations for form factors 
(\emph{cf.}~\cite[sect.~6]{Sm}).

\begin{rmk}
All $\tau$-functions satisfy bilinear Hirota equations (of KP
hierarchy) related to the time evolution $e^{H(t)}$. These equations
characterize the Sato Grassmannian itself.
\end{rmk} 

Some remarks can be made, however, in our infinite-dimensional setting. 
The ind-varieties $\mathcal{GR}$ and $\mathcal{GR}^{(n)}$ have 
two equivalent ind-structures given by pole orders and by Schubert 
varieties \cite{Ku}. Denote by $\mathcal{GR}=\varinjlim \mathcal{GR}_i$ 
one of these structures and let $\mathcal{D}(\mathcal{GR})=\varinjlim
\mathcal{D}(\mathcal{GR}_i)$ be the sheaf of differential operators
on $\mathcal{GR}$. Then one considers right 
$\mathcal{D}=\mathcal{D}(\mathcal{GR})$-modules in the abelian category 
$\mathcal{M}(\mathcal{GR},\mathcal{O})$ of $\mathcal{O}^!$-modules 
(\emph{cf.}~\cite[sect.~7.11]{BD:qH} for a detailed exposition).

Quantum correlations functions constructed from intertwiners of 
$U_q(\widehat{\mathfrak{sl}}_n)$-modules should satisfy quantum 
differential equations in the sense of Frenkel-Reshetikhin
\cite{FR}. We use the notation $\mathcal{D}_q(\tau_W^{})$ 
for associated \emph{quantum $\mathcal{D}$-modules} on
$\mathcal{GR}^{(n)}$ parametrized by affine local $\GL_n$-opers.

\section{Affine Miura opers and Miura transformation}

We use \cite[7.2 and 8.3.5]{BZF} as a reference for affine Miura opers.
Let $G$ be a reductive group, $X$ a smooth complex curve and $P\in|X|$.
Denote $\mathrm{L}G_X^-\subset \mathrm{L}G$ the subgroup consisting of 
loops which extends to $X\setminus P$.

Let $B$ be a Borel subgroup of $G$ and $B^-$ a transverse Borel subgroup
(such that $[B^-]\in G/B$ lies in the open $B$-orbit). Choose a point $Q$
on $X$ distinct from $P$ and consider the subgroup 
$\mathrm{L}G_{\,X}^{\scriptscriptstyle{--}}\subset 
\mathrm{L}G_X^-$ consisting of ``negative'' loops 
whose values at $Q$ lie in $B^-$.

\begin{dfn}
A \emph{local affine Miura $G$-oper} is a quadruple $(\mathcal{U},
\nabla,\mathcal{U}^{\scriptscriptstyle{--}},\mathcal{U}^+)$ where 
$\mathcal{U}$ is an $\mathrm{L}G$-torsor on $\widehat{D}_t$ with a 
connection $\nabla$, a flat reduction $\mathcal{U}^{\scriptscriptstyle{--}}$ 
to $\mathrm{L}G_{\,X}^{\scriptscriptstyle{--}}$ and a reduction 
$\mathcal{U}^+$ to $\mathrm{L}G^+$ in ``tautological relative position'' 
with respect to $\nabla$ (\emph{loc.~cit.}). 
\end{dfn}

\begin{rmk} One can define local affine $G$-opers in the same vein replacing
$\mathrm{L}G_{\,X}^{\scriptscriptstyle{--}}$ by $\mathrm{L}G_X^-$.
\end{rmk}

Local Miura $\GL_n$-oper is a connection of the type
\begin{equation}
\partial_t-
\begin{pmatrix}
\chi_1^{} &0 &0 &\cdots &0 \cr
1 &\chi_2^{} &0 &\cdots &0 \cr
0 &1 &\ddots &\ddots &\vdots \cr
\vdots &\ddots &\ddots &\ddots &0 \cr
0 &\cdots &0 &1 &\chi_n^{}
\end{pmatrix},
\; \chi_i^{}\in\mathbb{C}[\![t]\!],
\end{equation}
on a vector bundle of rank $n$ over $\widehat{D}_t$. The \emph{Miura
transformation} is the transformation from this connection to the unique
gauge equivalent connection of the form \eqref{eq:locoper}. Geometrically,
it may be described as follows \cite[3.8]{Fr:aff}. Let $(\mathcal{E},\nabla)$
be a local $\GL_n$-oper with another full flag of subbundles 
\begin{equation}
\mathcal{E}_1'\subset \mathcal{E}_2'\cdots\subset\mathcal{E}_n'
\end{equation} 
preserved by $\nabla$. Such a flag defines a connection $\nabla'$ on 
$\gr(\mathcal{E})=\oplus_{i=1}^n \gr_i^{}(\mathcal{E})$. Then the Miura 
transformation associates $(\mathcal{E},\nabla)$ to
$\big(\gr(\mathcal{E}'),\nabla'\big)$. 

\medskip Now let $\big(\widehat{\mathcal{E}},\widehat{\nabla}\big)$
be a local affine $\GL_n$-oper and
\begin{equation}
\widehat{\mathcal{E}}=\widehat{\mathcal{E}}_0'\subset
\widehat{\mathcal{E}}_1'\subset\widehat{\mathcal{E}}_2'\subset\cdots
\subset\widehat{\mathcal{E}}_n'=\widehat{\mathcal{E}}'(Q)
\end{equation}
another parabolic structure preserved by $\widehat{\nabla}$. 
As above, the \emph{affine Miura transformation} associates 
$\big(\widehat{\mathcal{E}},\widehat{\nabla}\big)$ to 
$\big(\gr\big(\widehat{\mathcal{E}'}\big),\widehat{\nabla}'\big)$.

\begin{prop}
Let $\big(\gr\big(\widehat{\mathcal{E}'}\big),\widehat{\nabla}'\big)$
be a local affine Miura $\GL_n$-oper and $\big(\widehat{\mathcal{E}},
\widehat{\nabla}\big)$ the associated local affine $\GL_n$-oper. Let
$W$ be a point of the $n$-reduced Sato Grassmannian $\mathcal{GR}^{(n)}$
which is the image of  $\big(\widehat{\mathcal{E}},\widehat{\nabla}\big)$
under the Krichever map. Then 
\begin{equation}
\big(\gr\big(\widehat{\mathcal{E}'}\big),\widehat{\nabla}'\big)
\mapsto \big(W,\widehat{\mathcal{E}}_1'/\widehat{\mathcal{E}}'
\subset \widehat{\mathcal{E}}_2'/\widehat{\mathcal{E}}'\subset
\cdots\subset \widehat{\mathcal{E}}_n'/\widehat{\mathcal{E}}')
\end{equation}
defines a point of the affine flag variety $\mathcal{FL}^{(n)}$.
\end{prop}

\mn \emph{Proof.} Recall that there is a natural fibration
$\mathcal{FL}^{(n)}\rightarrow\mathcal{GR}^{(n)}$ whose
fiber is the usual flag variety $\GL_n/B$. Thus, the pair
\begin{equation}
\big(W,\widehat{\mathcal{E}}_{\bullet}'/\widehat{\mathcal{E}}'\big)=
\big(W,\widehat{\mathcal{E}}_1'/\widehat{\mathcal{E}}'
\subset \widehat{\mathcal{E}}_2'/\widehat{\mathcal{E}}'\subset
\cdots\subset \widehat{\mathcal{E}}_n'/\widehat{\mathcal{E}}'\big)
\end{equation}
defines, indeed, a point of $\mathcal{FL}^{(n)}$ under this fibration.

\section{$\mathcal{D}$-modules for affine Toda theories}

Consider the generating functions 
\begin{equation}
\psi^+(\mathfrak z)=\sum_{i\in\mathbb{Z}} \psi_i^+ \mathfrak{z}^i
\textrm{ and }
\psi^-(\mathfrak z)=\sum_{i\in\mathbb{Z}} \psi_i^- \mathfrak{z}^{-i}
\end{equation}
of free fermions $\psi_i^{\pm}$ indexed by $i\in\mathbb{Z}$ 
(\emph{cf.}~sect.~\ref{sect:qwedges}). Introduce ``positive'' and ``negative" 
times $t=(t_1,t_2,\cdots)$ and $t'=(t_1',t_2',\cdots)$ resp., and
denote $\xi(t,\mathfrak{z})=\sum_{n>0} t_n\mathfrak{z}^n$. The 
\emph{two-dimensional Toda theory} is related to the following
time evolution $g(t,t')$:
\begin{align}
\psi^+(\mathfrak{z}) &\mapsto\mathfrak{z}^n e^{\xi(t,\mathfrak{z})+
\xi(t',\mathfrak{z}^{-1})} \psi(\mathfrak{z})\\
\psi^-(\mathfrak{z}) &\mapsto\mathfrak{z}^{-n} e^{-\xi(t,\mathfrak{z})
-\xi(t',\mathfrak{z}^{-1})} \psi(\mathfrak{z})
\end{align}
(with essential singularities at $\mathfrak{z}=0,\infty$) of
$g=\exp\big(\sum_i a_i^{}\psi^+(p_i^{})\psi^-(q_i^{})\big)$
\cite[\S9]{JM}. The \emph{$\tau$-functions} are defined as before by
\begin{equation}
\tau_g^{}(t,t')=\langle\,0\,|\,g(t,t')\,|\,0\,\rangle
\end{equation}
(loc.~cit.). The usual reduction to $\textrm{A}_n^{(1)}$ \cite[\S8,9]{JM}
defines $\tau$-functions of affine Toda chains. 

\medskip Let $\big(\gr\big(\widehat{\mathcal{E}'}\big),
\widehat{\nabla}'\big)$ be a local affine Miura $\GL_n$-oper
and denote $\big(W,\widehat{\mathcal{E}}_{\bullet}'/
\widehat{\mathcal{E}}'\big)$ the associated point of 
$\mathcal{FL}^{(n)}$. For any $\tau$-function 
\begin{equation}
\tau_W^{}(t,t')=\langle\,0\,|\,g_W^{}(t,t')\,|\,0\,\rangle
\end{equation}
corresponding to $\big(W,\widehat{\mathcal{E}}_{\bullet}'/
\widehat{\mathcal{E}}'\big)$, we construct $\mathcal{D}$-modules
$\mathcal{D}(\tau_W^{})$ and $\mathcal{D}_q(\tau_W^{})$ as in
sect.~\ref{sect:qwedges}.

\section{Proof of the main theorem}

The result follows easily from the constructions above.

\medskip Let $\big(\gr\big(\widehat{\mathcal{E}'}\big),
\widehat{\nabla}'\big)$ be a local affine Miura $\GL_n$-oper
and $\big(W,\widehat{\mathcal{E}}_{\bullet}'/
\widehat{\mathcal{E}}'\big)$ the associated point of 
$\mathcal{FL}^{(n)}$. Under the Miura transformation,
we get a local affine $\GL_n$-oper $(\mathcal{E},\nabla)$.
It is easy to see that the Miura transform of
$\big(W,\widehat{\mathcal{E}}_{\bullet}'/
\widehat{\mathcal{E}}'\big)$ coincides with the
point $W\in\mathcal{GR}^{(n)}$ obtained from
$(\mathcal{E},\nabla)$ via the Krichever map.

\medskip Moreover, the $\tau$-function $\tau_W^{}(t,t')$
of the affine Toda $\widehat{\mathfrak{sl}}_n$-chain becomes the 
$\tau$-function $\tau_W^{}(t)$ of the $n$th KdV hierarchy. 
Clearly, this transformation preserves null vectors. Thus,
the Givental type $\mathcal D$-modules $\mathcal{D}(\tau_W^{})$ and 
$\mathcal{D}_q(\tau_W^{})$ on $\mathcal{FL}^{(n)}$ are
transformed to Beilinson-Drinfeld type  $\mathcal D$-modules
on $\mathcal{GR}^{(n)}$.

\section{Alternative construction of $\mathcal{D}$-modules}

In this section we briefly indicate how one can construct $\mathcal{D}$-modules
on $\mathcal{GR}^{(n)}$ and $\mathcal{FL}^{(n)}$ using a local affine
version of the Beilinson-Drinfeld approach.

\smallskip In the case of the $n$th KdV hierarchy, an analogue of the 
Hitchin fibration is the fibration $S\rightarrow\mathcal{AO}^{(n)}$ whose
fibers are (Jacobians) of spectral curves. Here $\mathcal{AO}^{(n)}$ denotes
the moduli space of local affine $\GL_n$-opers (so-called \emph{abelianized
Grassmannian}). Recall (\emph{cf.}~\cite{AB,LM}) that any local affine 
$\GL_n$-oper $m\in\mathcal{M}$, \emph{i.~e.}, any Krichever module, defines 
a spectral curve $X_m$. The fibration above is closely related to 
$T^*\mathcal{GR}^{(n)}$ and one can construct Beilinson-Drinfeld type 
$\mathcal{D}$-modules on $\mathcal{GR}^{(n)}$ indexed by local affine 
$\GL_n$-opers.

\medskip In the case of the affine Toda $\widehat{\mathfrak{sl}}_n$-chain,
an analogue of the Hitchin fibration is the fibration 
$S\rightarrow\mathcal{MO}^{(n)}$ whose fibers are (Jacobians)
of spectral curves of this chain \cite[sect.~1]{McD}. Here 
$\mathcal{MO}^{(n)}$ denotes the moduli space of local affine 
Miura $\GL_n$-opers. Using the Beilinson-Drinfeld procedure one can 
construct $\mathcal{D}$-modules on $\mathcal{FL}^{(n)}$ indexed by 
such opers.

\bn {\footnotesize \textbf{Acknowledgements.} I am grateful to Mark
Spivakovsky for useful remarks.} 

\nocite{*}
\bibliographystyle{amsplain}
\bibliography{qldms}

\bigskip
\begin{flushright}
Igor Potemine\\
Laboratoire Emile Picard\\
Universit\'e Paul Sabatier\\
118, route de Narbonne\\
31062 Toulouse C\'edex
\end{flushright}

\begin{flushright}
e-mail : potemine@picard.ups-tlse.fr
\end{flushright}

\end{document}